\documentclass[12pt]{article}
\usepackage[latin1]{inputenc}
\usepackage{amsmath}
\usepackage{amssymb}
\usepackage{esint}

\makeatletter
\usepackage{latexsym}
\usepackage{enumerate}
\usepackage{epsf}
\usepackage{graphicx}
\usepackage{comment}
\usepackage{appendix}

\newfont{\bb}{msbm10 at 12pt}
\newfont{\tbb}{msbm10 at 8pt}

\def\r{\hbox{\bb R}}

\def\h{\hbox{\bb H}}
\def\n{\hbox{\bb N}}

\def\c{\hbox{\bb C}}
\def\s{\hbox{\bb S}}
\def\tl{\hbox{\tbb L}}
\def\ts{\hbox{\tbb S}}

\def\tr{\hbox{\tbb R}}

\def\R{\r}
\def\H{\mathbb{H}}

\def\stop{\hfill$\Box$}

\newcommand{\set}[1]{\left\{#1\right\}}
\newcommand{\meta}[2]{\langle #1,#2 \rangle }

\usepackage{amsthm}
\usepackage{times}
\usepackage{amscd}
\usepackage{epsf}

\numberwithin{equation} {section}

\makeatother

\begin{document}

\theoremstyle{plain}\newtheorem{lemma}{Lemma}[section] 
\theoremstyle{plain}\newtheorem{proposition}{Proposition}[section]
\theoremstyle{plain}\newtheorem{theorem}{Theorem}[section] 
\theoremstyle{plain}\newtheorem*{main theorem}{Main Theorem} 
\theoremstyle{plain}\newtheorem{example}{Example}[section]
\theoremstyle{plain}\newtheorem{remark}{Remark}[section] 
\theoremstyle{plain}\newtheorem{corollary}{Corollary}[section]
\theoremstyle{plain}\newtheorem{definition}{Definition}[section]
\theoremstyle{plain}\newtheorem{acknowledge}{Acknowledgment}
\theoremstyle{plain}\newtheorem{conjecture}{Conjecture}
\begin{center}
\rule{15cm}{1.5pt} \vspace{0.4cm}

\textbf{\Large{}{}On nonnegatively curved hypersurfaces in $\H^{n+1}$}{\Large{} }

\vspace{0.4cm}

{\large{}{}Vincent Bonini$\,^{\dag}$, Shiguang Ma$\,^{\dag}$}\footnote{{Corresponding author. The author is support by NSFC grant No.11301284
and NSFG grant No.11571185.}}{\large{}{}, and Jie Qing$\,^{\star}$}\footnote{{The author is partially supported by NSF DMS-1303543}}\\

\vspace{0.3cm}
 \rule{15cm}{1.5pt} 
\par\end{center}

\vspace{0.2cm}
\noindent$\mbox{}^{\dag}$ Department of Mathematics, Cal Poly State University,
San Luis Obispo, CA 93407; \\
e-mail: vbonini@calpoly.edu \vspace{0.2cm}

\noindent $\mbox{}^{\ddag}$ Department of Mathematics and LPMC, Nankai
University, Tianjin, China; \\
 e-mail: msgdyx8741@nankai.edu.cn \vspace{0.2cm}

\noindent $\mbox{}^{\star}$ Department of Mathematics, University
of California, Santa Cruz, CA 95064; \\
 e-mail: qing@ucsc.edu 
\begin{abstract}
In this paper we prove a conjecture of Alexander and Currier that
states, except for covering maps of equidistant surfaces in hyperbolic
3-space, a complete, nonnegatively curved immersed hypersurface in
hyperbolic space is necessarily properly embedded. 
\end{abstract}

\section{Introduction}

\label{Sect:Intro}

Suppose that $\phi:M^{n}\to\mathbb{R}^{n+1}$ is an immersed hypersurface
with principal curvatures $\kappa_{1},\dots,\kappa_{n}$. Then $\phi$
is said to be 
\begin{itemize}
\item convex at a point if $\kappa_{i}\geq0$ for all $i=1,\dots,n$. 
\item of nonnegative Ricci curvature if $\kappa_{i}(\sum_{k=1}^{n}\kappa_{k})\geq\kappa_{i}^{2}$
for all $i=1,\dots,n$. 
\item nonnegatively curved if $\kappa_{i}\kappa_{j}\geq0$ for all $i,j=1,\dots,n$. 
\end{itemize}
It is easily seen that up to orientation all three of the curvature
conditions above are pointwise equivalent for hypersurfaces immersed in
Euclidean space. An immersed hypersurface in Euclidean space is said
to be locally convex if the hypersurface is locally supported by a
hyperplane. It is not true that nonnegativity of the sectional curvatures
alone implies local convexity of a hypersurface (cf. \cite{Sack}).
\\

The study of nonnegatively curved immersed hypersurfaces goes back
to Hadamard, who showed that a compact, strictly convex, immersed
surface in Euclidean 3-space is necessarily embedded \cite{Hdm}.
This result was later extended in \cite{Stoker,VH,cLashof,Sack} to such
that a complete, nonnegatively curved, nonflat, immersed hypersurface
in Euclidean space is necessarily embedded as a boundary of convex
body. \\

In this paper we consider oriented immersed hypersurfaces $\phi:M^{n}\to\mathbb{H}^{n+1}$
in hyperbolic space. The following pointwise curvature conditions
are no longer equivalent: 
\begin{itemize}
\item (strictly) convex at a point if $\kappa_{i} > 0$ for all $i=1,\dots,n$. 
\item nonnegative Ricci curvature if $\kappa_{i}(\sum_{k=1}^{n}\kappa_{k})\geq n-1+\kappa_{i}^{2}$
for all $i=1,\dots,n$. 
\item nonnegatively curved if $\kappa_{i}\kappa_{j}\geq1$ for all $i,j=1,\dots,n$. 
\item (non-strictly) horospherically convex if $\kappa_{i}\geq1$ for all
$i=1,\dots,n$. 
\end{itemize}
In fact, they are in strictly ascending order as listed above (cf.
\cite{Eps1,Eps2,AlCu,AlCu2}). Do Carmo and Warner \cite{dCW} showed that a
compact, convex, immersed hypersurface in hyperbolic space is necessarily
embedded. For noncompact cases, even with strict convexity, a complete,
immersed hypersurface in hyperbolic space may not be embedded \cite{Eps2}
(see also \cite{Spv}, pg. 84). On the other hand, Currier \cite{Cu}
showed that a (non-strictly) horospherically convex, complete, immersed
hypersurface in hyperbolic space is necessarily embedded and, if noncompact,
a horosphere. Therefore one wonders whether a complete immersed hypersurface
with nonnegative sectional curvature or even nonnegative Ricci curvature is
necessarily embedded? \\

Naturally the embeddedness problem for a complete noncompact hypersurface
in hyperbolic space is related to its asymptotic boundary at infinity.
The asymptotic boundaries at infinity of complete hypersurfaces with
nonnegative curvature in hyperbolic space have been studied in \cite{Eps2,AlCu,AlCu2}.
In \cite{Eps2}, using hyperbolic Gauss maps and the geometry of horospheres, Epstein showed that a complete embedding of $\mathbb{R}^2$ 
into $\mathbb{H}^3$ with nonnegative Gaussian curvature has a single point asymptotic boundary at infinity. Epstein also showed \cite{Eps2} 
that a complete, strictly convex, immersed surface in $\h^{3}$ with a single point asymptotic boundary at infinity is necessarily embedded 
as the analog of van Heijenoort's theorem \cite{VH} in hyperbolic 3-space. Epstein then asked if a complete immersed surface
in $\mathbb{H}^3$ with nonnegative Gaussian curvature is necessarily embedded. \cite{Eps2}. 
\\

Later in \cite{AlCu, AlCu2}, based on a theorem of Volkov and Vladimirova \cite{VV}
and the splitting theorem of Cheeger and Gromoll \cite{CG}, Alexander
and Currier showed that a complete, nonnegatively curved, embedded
hypersurface in hyperbolic space has an asymptotic boundary at infinity
of at most two points and, if two points, it is an equidistant surface
about a geodesic line.  Alexander and Currier then in \cite{AlCu2} gave the 
precise statement of the conjecture as: Except for covering maps of equidistant surfaces in $\h^{3}$,
every nonnegatively curved immersed hypersurface in $\h^{n+1}$ is properly embedded. 
They also mentioned a sketch of a proof of this
conjecture for higher dimensions ($n\geq3$) suggested by Gromov.
Their conjecture remains completely open in the case when $n=2$. 
\\

In this paper we present proofs of the conjecture of Alexander and
Currier for the case when $n=2$ as well as all higher dimensions
($n\geq3$). Our main theorem is as follows:

\begin{main theorem}\label{Thm:MainThm} Except for covering maps
of equidistant surfaces in $\h^{3}$, a complete, nonnegatively curved,
immersed hypersurface in hyperbolic space $\h^{n+1}$ for $n\geq2$
is properly embedded. \end{main theorem}

Our approach for solving the conjecture of Alexander and Currier in
higher dimensions ($n\geq3$) is based on the recent work \cite{BEQ2}
on weakly horospherically convex hypersurfaces in hyperbolic space 
(cf. Definition \ref{horocon}), which may be considered as an extension 
of the embedding theorem in \cite{Eps2}. Please see Theorem 
\ref{embedding} and Theorem \ref{embedding-one} in Section \ref{Sect:LCT}. This approach 
was initiated by Epstein in \cite{Eps2}. One key issue is to derive the injectivity of the hyperbolic
Gauss map. We will rely on the injectivity theorem of Schoen and Yau \cite{SY,SY1}, while
Epstein \cite{Eps2} used the embeddedness. The other key issue is the size estimate for the asymptotic 
boundary at infinity. We will rely on the Hausdoff dimension estimate of Zhu \cite{Zhu}, while Epstein's approach in  \cite{Eps2} is based on similar results of Huber \cite{Hu} for subharmonic functions.  
\\

To prove the conjecture of Alexander and Currier in dimension $n=2$, we first establish a new proof 
of the classical result of Volkov and Vladimirova \cite{VV}, which states that the
only way to isometrically immerse the Euclidean plane $\R^{2}$ in $\h^{3}$ is as a covering map of an equidistant surface about a geodesic
line or as a horosphere. Our proof of the main theorem is then based on the sharp growth estimate \eqref{taliaferro} in Lemma \ref{Hu-Tal} 
for solutions to Gaussian curvature equations based on \cite{Hu,Tal, Tal-1}. 
The key lower bound estimate for solutions to Gaussian curvature equations, which is needed to use \cite{Tal, Tal-1}, is based on the non-collapsing result of Croke and Karcher \cite{Croke and Karcher} and a Harnack-type estimate from Li and Schoen \cite{Li and Schoen}. Our approach in spirit is to show that a complete, noncompact, nonnegatively curved, nonflat, immersed surface in $\h^{3}$ lies inside a horosphere, hence has an asymptotic boundary at infinity of exactly one point. Then the embeddedness follows from Epstein \cite{Eps2}. 
\\

This paper is organized as follows: In Section \ref{Sect:LCT} we
introduce the geometry of horospherical metrics for weakly horospherically
convex hypersurfaces in hyperbolic space and some framework from \cite{EGM,BEQ,BEQ2}.
In Section \ref{Sect:EmbHighDim} we apply the embedding Theorems
\ref{embedding} and \ref{embedding-one} (see also \cite{BEQ2})
to prove the conjecture of Alexander and Currier \cite{AlCu2} in
higher dimensions ($n\geq3$). In Section \ref{Sect:EmbSurfs} we
present the proof of the conjecture of Alexander and Currier \cite{AlCu2}
in the case when $n=2$. \\

\noindent \textbf{Acknowledgment} The authors would like to express
their gratitude to Professor Jose Espinar at IMPA for his interest
in this work. We are very appreciative of his careful reading that
led to the current version.


\section{Hyperbolic Gauss Maps and Horospherical Metrics}

\label{Sect:LCT}

In this section we recall the definitions of hyperbolic Gauss maps
and weak horospherical convexity to set our terminologies and notations.
Readers are referred to the papers \cite{Eps1,Eps3,EGM,BEQ,BEQ2}
for more details. \\

For $n\geq2$, we denote Minkowski spacetime by $\r^{1,n+1}$, which
is the vector space $\r^{n+2}$ endowed with the Minkowski spacetime
metric $\meta{}{}$ given by 
\[
\meta{\bar{x}}{\bar{x}}=-x_{0}^{2}+\sum_{i=1}^{n+1}x_{i}^{2},
\]
where $\bar{x}\equiv(x_{0},x_{1},\dots,x_{n+1})\in\r^{n+2}$. Then
hyperbolic space, de Sitter space, and the positive null cone are
given by the respective hyperquadrics 
\[
\begin{split}\h^{n+1} & =\set{\bar{x}\in\r^{1,n+1}:\,\meta{\bar{x}}{\bar{x}}=-1,\,x_{0}>0},\\
\s^{1,n} & =\set{\bar{x}\in\r^{1,n+1}:\,\meta{\bar{x}}{\bar{x}}=1},\\
\n_{+}^{n+1} & =\set{\bar{x}\in\r^{1,n+1}:\,\meta{\bar{x}}{\bar{x}}=0,\,x_{0}>0}.
\end{split}
\]
We identify the ideal boundary at infinity $\partial_{\infty}\h^{n+1}$
of hyperbolic space with the unit round sphere $\s^{n}$ sitting at
$x_{0}=1$.

\begin{definition}(cf. \cite{Br,Eps1,Eps3})\label{h-gauss} Let
$\phi:M^{n}\to\h^{n+1}$ denote an immersed oriented hypersurface
in $\h^{n+1}$ with unit normal $\eta:M^{n}\to\s^{1,n}$. The hyperbolic
Gauss map 
\[
G:M^{n}\to\s^{n}
\]
of $\phi$ is defined as follows: for $p\in M^{n}$, the image $G(p)\in\s^{n}$
is the point at infinity of the unique horosphere in $\h^{n+1}$ passing
through $\phi(p)$ and whose outward unit normal at $\phi(p)$ agrees
with $\eta(p)$. \end{definition}

Given an oriented, immersed hypersurface $\phi:M^{n}\to\h^{n+1}$
with unit normal vector field $\eta:M^{n}\to\s^{1,n}$, the light
cone map $\psi$ associated to $\phi$ is defined by 
\[
\psi:=\phi-\eta:M^{n}\to\n_{+}^{n+1}.
\]
As the ideal boundary $\s^{n}$ of $\h^{n+1}$ is identified with
the unit round sphere at $x_{0}=1$, we have 
\begin{equation}
\psi=e^{\rho}(1,G),\label{psiG}
\end{equation}
where $\psi_{0}=e^{\rho}$ is the so-called horospherical support
function of the hypersurface $\phi$ \cite{EGM}. Note that, in our
convention given in Definition \ref{h-gauss}, horospheres with outward
orientation are the unique surfaces such that both the hyperbolic
Gauss map and the associated light cone map are constant. Moreover,
if $x\in\s^{n}$ is the point at infinity of such a horosphere, then
$\psi=e^{\rho}(1,x)$ where $\rho$ is the signed hyperbolic distance
of the horosphere to the point $\mathcal{O}=(1,0,\dots,0)\in\h^{n+1}\subseteq\r^{1,n+1}$.
\\

Considering the fact that horospheres are intrinsically flat, one
can then use horospheres to define concavity/convexity for hypersurfaces
in hyperbolic space.

\begin{definition} (cf. \cite{Sc,EGM,BEQ2})\label{horocon} Let
$\phi:M^{n}\to\h^{n+1}$ be an immersed oriented hypersurface and
let $\mathcal{H}_{p}$ denote the horosphere in $\h^{n+1}$ that is
tangent to $\phi(M)$ at $\phi(p)$ whose outward unit normal at $\phi(p)$
agrees with $\eta(p)$. We will say that $\phi$ is weakly horospherically
convex at $p$ if there exists a neighborhood $V\subset M^{n}$ of
$p$ so that $\phi(V\backslash\{p\})$ stays outside of $\mathcal{H}_{p}$.
Moreover, the distance function of the hypersurface to the horosphere
does not vanish up to the second order at $p$ in any direction. \end{definition}

Due to \cite{EGM}, we have the following characterization of weakly
horospherically convex hypersurfaces.

\begin{lemma}(\cite{EGM})\label{hc} Let $\phi:M^{n}\to\h^{n+1}$
be an immersed oriented hypersurface. Then $\phi$ is weakly horospherically
convex at $p$ if and only if the principal curvatures $\kappa_{1},\dots,\kappa_{n}$
of $\phi$ at $p$ are simultaneously $>-1$. In particular, $\phi$
is weakly horospherically convex at $p$ implies that $dG$ is invertible
at $p$ and therefore the hyperbolic Gauss map of $\phi$ is a local
diffeomorphism. \end{lemma}

To realize this second statement, let $\{e_{1},\dots,e_{n}\}$ denote
an orthonormal basis of principal directions of $\phi$ at $p$ and
let $\kappa_{1},\dots,\kappa_{n}$ denote the associated principal
curvatures. Then $d\phi(e_{i})=e_{i}$ and $d\eta(e_{i})=-\kappa_{i}e_{i}$
for $i=1,\dots,n$, so as in \cite{EGM}, it follows that 
\begin{equation}
\meta{(d\psi)_{p}(e_{i})}{(d\psi)_{p}(e_{j})}_{\tr^{1,n+1}}=(1+\kappa_{i})^{2}\delta_{ij}=e^{2{\rho}}\meta{(dG)_{p}(e_{i})}{(dG)_{p}(e_{j})}_{\ts^{n}},\label{MetricRelation}
\end{equation}
where $g_{\ts^{n}}$ denotes the round metric on $\s^{n}$. Now given
an immersed oriented weakly horospherically convex hypersurface $\phi:M^{n}\to\h^{n+1}$,
one can use the hyperbolic Gauss map (or light cone map) to induce
a canonical locally conformally flat metric on $M^{n}$ as follows:

\begin{definition} (\cite{Eps2,Eps3,EGM}\label{Def:HoroMetric})
Let $\phi:M^{n}\to\h^{n+1}$ be an immersed oriented weakly horospherically
convex hypersurface. Then the hyperbolic Gauss map $G:M^{n}\to\s^{n}$
is a local diffeomorphism. We consider the locally conformally flat
metric 
\begin{equation}
g_{h}=\Psi^{*}\meta{}{}_{\tl^{n+2}}=e^{2\rho}G^{*}g_{\ts^{n}}\label{horom}
\end{equation}
on $M^{n}$ and call it the horospherical metric associated to the
immersed oriented weakly horospherically convex hypersurface $\phi$.
\end{definition}

For a weakly horospherically convex hypersurface $\phi$, its associated
light cone map $\Psi$ is spacelike and parameterizes a codimension
$2$ submanifold in $\r^{1,n+1}$. $\phi$ and $\eta$ provide two
unit normal fields to $\Psi$ and the second fundamental form is given
by 
\begin{equation}
II_{\Psi}(e_{i},e_{j})=(\frac{1}{1+\kappa_{i}}\phi-\frac{\kappa_{i}}{1+\kappa_{i}}\eta)g_{h}(e_{i},e_{j})\label{2FundForm}
\end{equation}
where $\{e_{1},\dots,e_{n}\}$ is an orthonormal basis of principal
directions with respect to $\phi$. Hence, due to the Gauss equations
in $\r^{1,n+1}$, the sectional curvatures of the horospherical metric
$g_{h}$ on $M^{n}$ are given by 
\begin{equation}
K_{g_{h}}(\frac{e_{i}}{1+\kappa_{i}},\frac{e_{j}}{1+\kappa_{j}})=1-\frac{1}{1+\kappa_{i}}-\frac{1}{1+\kappa_{j}}=\frac{\kappa_{i}\kappa_{j}-1}{(1+\kappa_{i})(1+\kappa_{j})}.\label{HorSectCurv}
\end{equation}
When $n\geq3$, the Schouten tensor then is given by 
\begin{equation}
Sch_{g_{h}}(e_{i},e_{j})=(\frac{1}{2}-\frac{1}{1+\kappa_{i}})g_{h}(e_{i},e_{j}).\label{HorSchouten}
\end{equation}
When $n=2$, instead, one considers the symmetric 2-tensor 
\begin{equation}
P=-\nabla_{G^{*}g_{\ts^{2}}}d\rho+d\rho\otimes d\rho-\frac{1}{2}(|d\rho|_{G^{*}g_{\ts^{2}}}^{2}-1)G^{*}g_{\mathbb{S}^{2}},
\end{equation}
whose eigenvalues are 
\begin{equation}
\frac{1}{2}-\frac{1}{1+\kappa_{1}}\ \text{ and }\ \frac{1}{2}-\frac{1}{1+\kappa_{2}},\label{p-eigen}
\end{equation}
whose trace is the Gaussian curvature 
\begin{equation}
K_{g_{h}}=\frac{\kappa_{1}\kappa_{2}-1}{(1+\kappa_{1})(1+\kappa_{2})},\label{p-trace}
\end{equation}
and whose divergence is $2dK_{g_{h}}$. Hence we get the Gaussian curvature
equation 
\begin{equation}
-\Delta_{G^{*}g_{\ts^{2}}}\rho+1=K_{g_{h}}e^{2\rho}.\label{gauss-curvature-equ}
\end{equation}

When the hyperbolic Gauss map $G:M^{n}\to\s^{n}$ of a weakly horospherically
convex hypersurface $\phi:M^{n}\to\h^{n+1}$ is injective, one may
push down the horospherical metric $g_{h}$ onto the image 
\begin{equation}
\Omega=G(M)\subset\s^{n}\label{gauss-image}
\end{equation}
to obtain the conformal metric 
\begin{equation}
\hat{g}_{h}=(G^{-1})^{*}g_{h}=e^{2\hat{\rho}}g_{\ts^{n}},\label{g-hat}
\end{equation}
where $\hat{\rho}=\rho\circ G^{-1}:\Omega\to\mathbb{R}$. When there
is no confusion, we will also refer to this conformal metric $\hat{g}_{h}$
as the horospherical metric. The correspondence between weakly horospherically
convex hypersurfaces $\phi:M^{n}\to\h^{n+1}$ in hyperbolic space
and the conformal metric $\hat{g}_{h}$ on the image $\Omega$ of
the Gauss map $G$ have been promoted in \cite{Eps2,EGM,BEQ,BEQ2, BQZ}.
The following result follows from the so-called global correspondence
from \cite{BEQ2,BQZ,EGM} and will be useful to our work here.

\begin{theorem}[cf. \cite{BEQ2, BQZ, EGM}] \label{Thm:GCT} For $n\geq2$,
let $\phi:M^{n}\to\h^{n+1}$ be a complete uniformly weakly horospherically
convex hypersurface with injective hyperbolic Gauss map $G:M^{n}\to\s^{n}$.
Then 
\begin{itemize}
\item $\phi$ induces a complete conformal metric $\hat{g}_{h}=e^{2\hat{\rho}}g_{\ts^{n}}$
on the image $\Omega=G(M)\subset\s^{n}$ with bounded curvature. 
\item More importantly, the asymptotic boundary $\partial_{\infty}\phi(M)\subset\s^{n}$
at infinity of the hypersurface $\phi$ in $\h^{n+1}$ coincides with
the boundary $\partial\Omega\subset\s^{n}$ of the Gauss map image. 
\item One may use the image $\Omega$ of Gauss map as the parameter space
to reparametrize $\phi$ so that the Gauss map 
\[
G(x)=x:\Omega\to\s^{n}
\]
and 
\begin{equation}
\phi_{t}=\frac{e^{\rho+t}}{2}(1+e^{-2(\rho+t)}(1+|\nabla\rho|^{2}))(1,x)+e^{-(\rho+t)}(0,-x+\nabla\rho)\label{normal-flow}
\end{equation}
is the nornal flow of the hypersurface $\phi(M)$. 
\end{itemize}
\end{theorem}

The contribution of \cite{BEQ} is the use of the normal flow of a
weakly horospherically convex hypersurface with injective hyperbolic
Gauss map to possibly unfold the hypersurface into an embedded one.
This is because the leaves of regular part of the normal flow are
the same as the level surfaces of the geodesic defining function of
the horospherical metric $\hat{g}_{h}$ (cf. \cite{BEQ,BEQ2}). For
instance, it is observed in \cite{BEQ} that any horospherical ovaloid
can be deformed along its normal flow into an embedded one. Consequently
this leads to new proofs of Obata type theorems for horospherical
ovaloids. In \cite{BEQ2,BQZ}, based on the global correspondence theorem,
we established an extension of the embedding theorem of Epstein \cite{Eps2}
as follows:

\begin{theorem}(cf. \cite{BEQ2, BQZ})\label{embedding} For $n\geq2$,
let $\phi:M^{n}\to\h^{n+1}$ be a complete uniformly weakly horospherically
convex hypersurface with injective hyperbolic Gauss map $G:M^{n}\to\s^{n}$.
Suppose that the asymptotic boundary $\partial_{\infty}\phi(M)$
at infinity of the hypersurface is a disjoint union of smooth closed
submanifolds in $\s^{n}$. Then, along the normal flow from the hypersurface,
the leaves eventually become embedded. \end{theorem}

An argument similar to those in \cite{VH,Eps2} results in the following
slight extension of the embedding theorem of Epstein \cite{Eps2}.

\begin{theorem} \label{embedding-one} For $n\geq2$, let $\phi:M^{n}\to\h^{n+1}$
be a complete, locally strictly convex, immersed hypersurface. Suppose
that the asymptotic boundary $\partial_{\infty}\phi(M)$ at infinity
of the hypersurface is a single point in $\s^{n}$. Then the hypersurface
is in fact embedded. \end{theorem}

\proof For convenience of readers, we would like to present a proof
based on the arguments in \cite{VH,Eps2}, which are similar to those
in \cite{BEQ2}. Since the asymptotic boundary at infinity of the
hypersurface is a single point in $\s^{n}$, one may find a family
of round $(n-1)$-spheres in $\s^{n}$ to foliate the sphere $\s^{n}$
with the point and its antipodal point deleted. Then the family of
hyperplanes whose asymptotic boundary at infinity are the family of
round $(n-1)$-spheres foliates hyperbolic space. To finish the argument
one simply needs to observe that, close to the first touch point of
the hyperplanes and the hypersurfaces from the antipodal point, the
hypersurface is locally embedded and the intersections of the hyperplanes
and hypersurfaces are embedded convex topological spheres. Moreover,
everything remains the same up to the end. The connectedness and convexity
of the hypersurface force each intersection to be connected and convex.
The embeddedness of the intersections is due to \cite{dCW}. \endproof


\section{Embeddedness in Higher Dimensions}

\label{Sect:EmbHighDim}

In this section we consider noncompact hypersurfaces immersed in hyperbolic
space with nonnegative sectional curvature and present a proof for
the conjecture of Alexander and Currier \cite{AlCu2} in higher dimensions
($n\geq3$). Based on the injectivity of development maps of Schoen
and Yau \cite{SY,SY1} and the Hausforff dimension estimates of Zhu
\cite{Zhu}, the proof of the conjecture of Alexander and Currier
\cite{AlCu2} is rather straightforward following our work in \cite{BEQ2}
and the brief summary in the previous section. \\

First of all, from the curvature relations \eqref{HorSectCurv}, we
have:

\begin{lemma} \label{curv-rel} Suppose that $\phi:M^{n}\to\h^{n+1}$
is a nonnegatively curved immersed hypersurface. Then $\phi$ is weakly
horospherically convex and the horospherical metric is also nonnegatively
curved. \end{lemma}

\proof It is easily seen that a nonnegatively curved hypersurface
in hyperbolic space is weakly horospherically convex, in fact, it
is strictly convex. Then the lemma is a simple consequence of \eqref{HorSectCurv}.
\endproof

There does not seem to be any analog of Lemma \ref{curv-rel} available
if we consider nonnegative Ricci curvature for the hypersurface $\phi$
instead. In higher dimensions ($n\geq3$), using the works in \cite{SY,SY1,Zhu},
we obtain the following:

\begin{lemma} \label{sy-zhu} For $n\geq3$, let $\phi:M^{n}\to\h^{n+1}$
be a complete, nonnegatively curved, immersed hypersurface. Then the
hyperbolic Gauss map is a development map from $(M^{n},\ g_{h})$
and injective. Moreover, the Hausdorff dimension of $\partial G(M)=\s^{n}\backslash G(M)$
is zero. 
\end{lemma}

\proof Due to the uniform weak horospherical
convexity (strict convexity) of the hypersurface $\phi$, the completeness
of the hypersurface implies the completeness of the horospherical
metric $g_{h}$. In the light of Lemma \ref{curv-rel}, $(M^{n},\ g_{h})$
is a complete, nonnegatively curved Riemannian manifold. Therefore
the lemma follows from the injectivity theorem of Schoen and Yau in
\cite{SY,SY1} and the Hausdorff dimension estimates of Zhu in \cite{Zhu}.
Notice that the theorem of Schoen and Yau only needs $g_{h}$ to have
nonnegative scalar curvature and the Hausdorff estimates of Zhu only
need $g_{h}$ to be Ricci nonnegative. 
\endproof

One more ingredient for our proof of the conjecture of Alexander and
Currier \cite{AlCu2} in higher dimensions ($n\geq3$) is the following:

\begin{lemma}\label{normal-flow-curv} Suppose that $\phi:M^{n}\to\h^{n+1}$
is a nonnegatively curved immersed hypersurface. Then along
the normal flow \eqref{normal-flow} the hypersurface remains nonnegatively
curved. \end{lemma}

\proof For the normal flow \eqref{normal-flow} in hyperbolic space,
one knows exactly how the principal curvatures evolve: 
\begin{equation}
\kappa_{i}^{t}=\frac{\kappa_{i}+\tanh{t}}{1+\kappa_{i}\tanh{t}}.\label{PrinCurvFlow}
\end{equation}
One may then calculate the sectional curvatures $K_{ij}^{t}=\kappa_{i}^{t}\kappa_{j}^{t}-1$
for $t>0$ to find 
\begin{equation}
K_{ij}^{t}=\kappa_{i}^{t}\kappa_{j}^{t}-1=\frac{K_{ij}(1-\tanh^{2}{t})}{(1+\kappa_{1}\tanh{t})(1+\kappa_{2}\tanh{t})}\geq0,\label{sect-curv-norm-flow}
\end{equation}
where $K_{ij}$ are the sectional curvatures of $\phi$.
\endproof

We are now ready to prove the conjecture of Alexander and Currier
\cite{AlCu2} in higher dimensions ($n\geq3$).

\vskip 0.1in \proof (\textit{Main Theorem in higher dimensions})
\quad{}For $n\geq3$, let $\phi:M^{n}\to\h^{n+1}$ be an immersed,
complete, noncompact hypersurface with nonnegative sectional curvature.
In the light of Lemma \ref{sy-zhu} the hyperbolic Gauss map $G:M^{n}\to\s^{n}$
is injective and the Hausdorff dimension of $\partial G(M)\subset\s^{n}$
is zero. According to Theorem \ref{Thm:GCT} (cf. \cite{BEQ2}), we
have 
\[
\partial_{\infty}\phi(M)=\partial G(M).
\]
Now, if $\partial_{\infty}\phi(M)=\partial G(M)$ were empty,
then $\phi(M)$ would be compact. Moreover, since any set of Hausdorff
dimension zero is totally disconnected, due to the splitting theorem
of Cheeger and Gromoll \cite{CG}, the asymptotic boundary $\partial_{\infty}\phi(M)=\partial G(M)$
consists of either one single point or exactly two points. \\

When $\partial_{\infty}\phi(M)$ is a single point, the result
follows from Theorem \ref{embedding-one}. Assume $\partial_{\infty}\phi(M)$
consists of exactly two points. We then first apply Theorem \ref{embedding}
(please also see \cite{BEQ2}) and observe that along the normal flow
the nonnegatively curved hypersurface $\phi_{t}$ is embedded for
sufficiently large $t$. Notice that the nonnegativity of the sectional
curvatures of $\phi_{t}$ follows from Lemma \ref{normal-flow-curv}.
Therefore, by the rigidity result of Alexander and Currier \cite{AlCu,AlCu2},
for $t$ sufficiently large the hypersurface $\phi_{t}$ has to be
an equidistant hypersurface about a geodesic line in hyperbolic space.
This forces the hypersurface $\phi$ to be an equidistant hypersurface
in hyperbolic space. Thus the proof of the conjecture of Alexander
and Currier \cite{AlCu2} in higher dimensions ($n\geq3$) is complete.
\endproof


\section{Embeddedness of Nonnegatively Curved Surfaces}

\label{Sect:EmbSurfs}

In this final section we consider noncompact, complete surfaces immersed in
$\h^{3}$ with nonnegative Gaussian curvature and present a proof
of the conjecture of Alexander and Currier \cite{AlCu2} in dimension
$2$. \\

Suppose that $\phi:M^{2}\to\h^{3}$ is a complete, nonnegatively curved,
immersed surface. We may assume the surface is locally strictly convex
after a change of orientation, if necessary. Therefore the hyperbolic
Gauss map $G:M^{2}\to\s^{2}$ is a local diffeomorphism, and the horospherical
metric $g_{h}$ is complete (cf. Theorem \ref{Thm:GCT}) and nonnegatively
curved in the light of \eqref{p-trace}. In fact, the symmetric tensor
$P$ associated with the horosphericl metric $g_{h}$ satisfies 
\begin{equation}
-\frac{1}{2}g_{h}<P<\frac{1}{2}g_{h}\label{curv-bound}
\end{equation}
according to \eqref{p-eigen}. With the complex structure given by
the horospherical metric $g_{h}$ the Gauss map $G$ is a conformal
map into the Riemann sphere. Lemma \ref{sy-zhu} breaks down in dimension
2 because of the abundance of local holomorphic functions (the lack
of Liouville Theorem). The search for a type of Picard theorem for
holomorphic functions analogous to Lemma \ref{sy-zhu} in dimension
2 is technically much more difficult, though it seems to be a classic
topic. We are going to rely on the growth estimate \eqref{taliaferro}
in Lemma \ref{Hu-Tal} based on \cite{Hu,Tal, Tal-1} for the support function
$\rho$ as a solution to the Gaussian curvature equation \eqref{gauss-curvature-equ}.
The novelty of our approach is to recognize that nonflatness implies
that the asymptotic boundary at infinity consists of exactly one point
and embeddedness then follows directly from the embedding theorem
of Epstein \cite{Eps2} as a hyperbolic analog of the embedding theorem
of van Heijenoort \cite{VH}. \\

Let $\pi:\widetilde{M}^{2}\to M^{2}$ be the universal covering map.
Then we consider the new parametrization $\tilde{\phi}=\phi\circ\pi:\widetilde{M}^{2}\to\h^{3}$
with the hyperbolic Gauss map $\tilde{G}=G\circ\pi:\widetilde{M}^{2}\to\s^{2}$
and the horospherical metric $\tilde{g}_{h}=\pi^{*}g_{h}$ whose Gaussian
curvature $K_{\tilde{g}_{h}}=K_{g_{h}}\circ\pi\geq0$. Most importantly
we have the symmetric tensor 
\[
\tilde{P}=P\circ\pi=-\nabla^2_{\tilde{G}^{*}g_{\ts^{2}}}\tilde{\rho}+d\tilde{\rho}\otimes d\tilde{\rho}-\frac{1}{2}(|d\tilde{\rho}|^{2}_{\tilde{G}^{*}g_{\ts^{2}}}-1)\tilde{G}^{*}g_{\ts^{2}},
\]
where $\tilde{\rho}=\rho\circ\pi$ and 
\begin{equation}
-\frac{1}{2}\tilde{g}_{h}<\tilde{P}<\frac{1}{2}\tilde{g}_{h}.\label{tensor-p-tilde}
\end{equation}
It follows from Theorem 15 in \cite{Hu} of Huber that $(M^{2},g_{h})$
is parabolic when the surface $\phi$ is nonnegatively curved. Therefore
the universal cover $\widetilde{M}^{2}$ of $M^{2}$ is biholomorphic
to the complex plane $\c$.


\subsection{Flat Cases}

In this subsection we present a proof to the following theorem of
Volkov and Vladimirova \cite{VV}. Our proof paves a way for us to handle the nonflat
cases in next subsection.

\begin{theorem}[{\cite{VV}}] Let $\phi$ be an isometric immersion
from Euclidean plane to hyperbolic 3-space. Then $\phi$ is either
a covering map of an equidistant surface about a geodesic line in
$\h^{3}$ or it is an embedded horosphere. \end{theorem}

\proof First of all it follows from \eqref{p-trace} that $K_{g_{h}}\equiv0$
whenever $K_{\phi}\equiv0$. Therefore $(\mathbb{R}^{2},\ g_{h})$
is isometric to the Euclidean plane. Let $z=(x,y)$ be the Euclidean
coordinate for $(\mathbb{R}^{2},\ g_{h})$ so that 
$$
|dz|^{2}= g_h = e^{2\rho}G^*g_{\mathbb{S}^2}.
$$
From the properties of the tensor $P$, we know that $P$ is a symmetric
2-tensor, which is trace-free, divergence-free and bounded in the
sense that 
\[
-\frac{1}{2}|dz|^{2}<P<\frac{1}{2}|dz|^{2}.
\]
Thus $P$ is in fact constant since $P$ is associated with a bounded holomorphic 
function on $\c$. \\

This implies that the principal curvatures of the surface are both
constant (i.e. the surface is an isoparametric surface). Therefore it is a horosphere when $P=0$
and an equidistance surface when $P\neq 0$ according to the classification of isoparametric 
surfaces in hyperbolic 3-space (cf. for example, \cite{Cartan,Cecil,VV}). So the proof is complete.
\endproof

We would like to point out that in the case when $P=0$ (i.e. when the surface is a horosphere), 
one in fact can explicitly find that
\begin{equation}
\rho(x,y)=\log(C[(x-x_{0})^{2}+(y-y_{0})^{2}]+\frac{1}{4C})\label{FlatHorSupp1}
\end{equation}
for some positive constant $C$. 


\subsection{Nonflat Cases}

In this subsection we consider a complete, noncompact, nonnegatively
curved, nonflat, immersed surface $\phi:M^{2}\to\h^{3}$. We will
focus on how to recognize and use the nonflatness. From Huber's result \cite{Hu},
we know the universal cover $(\widetilde{M}^{2},\ \tilde{g}_{h})$
is globally conformal to the Euclidean plane. Let $z=(x,y)$
be the Euclidean coordinate for $\widetilde{M}^{2}$ so that 
\[
e^{2\tilde{\rho}_{0}}|dz|^{2}=\tilde{g}_{h}= e^{2\tilde{\rho}}\tilde{G}^{*}g_{\ts^{2}}.
\]
Rewrite the relation above as 
\[
|dz|^{2}=e^{2(\tilde{\rho}-\tilde{\rho}_{0})}\tilde{G}^{*}g_{\ts^{2}}=e^{2\rho_{0}}\tilde{G}^{*}g_{\ts^{2}}
\]
for $\rho_0 = \tilde\rho - \tilde{\rho}_0$ and consider the symmetric 2-tensor 
\begin{equation}
P_{0}=-\nabla_{\tilde{G}^{*}g_{\ts^{2}}}^2\rho_{0}+d\rho_{0}\otimes d\rho_{0}-\frac{1}{2}(|d\rho_{0}|_{\tilde{G}^{*}g_{\ts^{2}}}^{2}-1)\tilde{G}^{*}g_{\ts^{2}}.\label{tensor-p-0}
\end{equation}
It is perhaps helpful to think that with the Gauss map $\tilde{G}$
and support function $e^{\rho_{0}}$, $P_{0}$ corresponds to a ``surface"
in $\h^{3}$ as in Theorem \ref{Thm:GCT}. What is this ``surface"? From
the discussion in the flat cases in the previous subsection we know that it is a horosphere if
$P_{0}$ vanishes. The following is a simple calculation.

\begin{lemma}\label{Bij-Pij relation} In the $(x,y)$ coordinates
\begin{align*}
(P_{0})_{11} & =\partial_{x}^{2}\tilde{\rho}_{0}-\frac{1}{2}((\partial_{x}\tilde{\rho}_{0})^{2}-(\partial_{y}\tilde{\rho}_{0})^{2})+{\tilde{P}}_{11},\\
(P_{0})_{22} & =\partial_{y}^{2}\tilde{\rho}_{0}-\frac{1}{2}((\partial_{y}\tilde{\rho}_{0})^{2}-(\partial_{x}\tilde{\rho}_{0})^{2})+{\tilde{P}}_{22},\\
(P_{0})_{12} & =(P_{0})_{21}=\partial_{x}\partial_{y}\tilde{\rho}_{0}-(\partial_{x}\tilde{\rho}_{0})(\partial_{y}\tilde{\rho}_{0})+{\tilde{P}}_{12},
\end{align*}
where 
\[
{\tilde{P}}=-\nabla_{{\tilde{G}}^{*}g_{\ts^{2}}}^{2}{\tilde{\rho}}+d{\tilde{\rho}}\otimes d{\tilde{\rho}}-\frac{1}{2}(| d\tilde{\rho}|_{{\tilde{G}}^{*}g_{\ts^{2}}}^{2}-1){\tilde{G}}^{*}g_{\ts^{2}}
\]
is the Schouten tensor for the surface $\tilde{\phi}$. \end{lemma}

The most important technical tool in this case is the following sharp
growth estimates for solutions to Gaussian curvature equations based
on \cite[Theorem 10]{Hu}, \cite[Lemma 3]{Tal}, and \cite[Theorem 2.1]{Tal-1}. We will present
the proof in the next subsection.

\begin{lemma}\label{Hu-Tal} Suppose that $(\mathbb{R}^{2},e^{2u}|dz|^{2})$
is complete, noncompact, nonnegatively curved, and nonflat. If the
Gaussian curvature is bounded, then 
\begin{equation}
u=-m\log\sqrt{1+|z|^{2}}+o(\log\sqrt{1+|z|^{2}})\quad\text{as \ensuremath{|z|\to\infty}}\label{taliaferro}
\end{equation}
for some $m\in(0,1]$. 

\end{lemma}

We are now ready to prove that $P_{0}$ vanishes.

\begin{lemma}\label{key} The Schouten tensor $P_{0}$ in \eqref{tensor-p-0}
vanishes identically on $\r^{2}$ and $\rho_{0}$ is given as a solution
in \eqref{FlatHorSupp1}. \end{lemma}

\proof First of all we know that $P_{0}$ is trace-free and divergence-free since $e^{2\rho_0}\tilde 
G^*g_{\mathbb{S}^2}=|dz|^2$ is flat. To show $P_{0}$ is in fact identically zero one just needs to show $|P_{0}|\in L^{p}(\mathbb{R}^{2})$
for some $p>1$, in the light of, for instance, \cite[Theorem 3]{Yau}. As the Gaussian
curvature of $\tilde{g}_{h}$ is bounded, by applying Lemma \ref{Hu-Tal},
we get 
\[
\tilde{\rho}_{0}=-m\log\sqrt{1+|z|^{2}}+o(\log\sqrt{1+|z|^{2}})\quad\text{as \ensuremath{|z|\to\infty}}
\]
for some $m\in(0,1]$. Then from \eqref{tensor-p-tilde} we know that
\begin{equation}
|\tilde{P}|\leq Ce^{2\tilde{\rho}_{0}}\leq\frac{C}{(1+|z|^{2})^{\frac{m}{2}}},\label{P-tilde estimate}
\end{equation}
and hence $|\tilde{P}|$ belongs to $L^{p}(\mathbb{R}^{2})$ for some
large $p>1.$ From interior estimates to the Gaussian curvature equation
\begin{equation}\label{Gauss-Equation}
-\Delta \tilde\rho_0 = K_{\tilde g_h} e^{2\tilde\rho_0}
\end{equation}
From the Schauder and $L^p$ estimates of \cite{Gt}, we have
\begin{equation}
\begin{cases}
R^{2-\frac{2}{p}}\|\partial^{2}\tilde{\rho}_{0}\|_{L^{p}(B_{R}(0))} & \leq 
C(\|\tilde{\rho}_{0}\|_{C^{0}(B_{2R}(0))}+R^{2-\frac{2}{p}}\|K_{\tilde g_h}e^{2\tilde{\rho}_{0}}\|_{L^{p}(B_{2R}(0))}),\\
r\|\partial\tilde{\rho}_{0}\|_{C^{0}(B_{r}(z))} & \leq C(\|\tilde{\rho}_{0}\|_{C^{0}(B_{2r}(z))}+
r^{2}\|K_{\tilde g_h}e^{2\tilde{\rho}_{0}}\|_{C^{0}(B_{2r}(z))}).
\end{cases}\label{derivatives estimate}
\end{equation}
From (\ref{P-tilde estimate}) and the first inequality of (\ref{derivatives estimate}) as $R\to\infty$,
we have $\partial^{2}\tilde{\rho}_{0}\in L^{p}(\mathbb{R}^{2})$ for
any $p$ sufficiently large since $K_{\tilde g_h}$ is bounded. Meanwhile, from the second inequality of 
\eqref{derivatives estimate}
and $m\in (0, 1]$ for 
$$r=(1+|z|^{2})^{\frac{m}{4}} < \frac 12 |z|,
$$ 
at least when $|z| > 2\sqrt{2}$,  we get 
$$|\partial\tilde{\rho}_{0}(z)|\leq\frac{C}{(1+|z|^{2})^{\frac{m}{4}}}(\log|z|+C),
$$
which implies that $|\partial\tilde{\rho}_{0}(z)|^{2}\in L^{p}(\mathbb{R}^{2})$
for $p$ sufficiently large. Therefore, due to Lemma 4.1, it follows that  $|P_0| \in L^p(\mathbb{R}^2)$.
\\

With $P_0=0$, we now know that the support function $\rho_0$ and the Gauss map $\tilde G$ indeed induce a real ``surface", which in fact is a horosphere. Thus the proof is complete.
\endproof

We are now ready to complete the proof of the conjecture of Alexander
and Currier \cite{AlCu2} in dimension 2.
\\

\noindent{\it Proof of the Main Theorem in nonflat cases in dimension 2} \quad
From Lemma \ref{key} we know $\tilde{G}$ is an injective map which
misses only one point $q\in\mathbb{S}^{2}.$ So the covering map $\pi$
is a diffeomorphism. From \eqref{FlatHorSupp1} and \eqref{taliaferro} we have
\begin{align*}
\tilde{\rho} & =\rho_{0}+\tilde{\rho}_{0}=(2-m)\log|z|+o(\log|z|),m\in(0,1]
\end{align*}
 as $z\rightarrow+\infty$. So $\tilde{\rho}(\tilde{G}^{-1}(\xi))\rightarrow+\infty$
as $\xi\rightarrow q$, which, together with the proof of Lemma 3.2
of \cite{BEQ2}, implies that $\partial_{\infty}\phi(M)=\{q\}$. We remark that one may 
derive the same conclusion from \cite{BQZ}. From
the embedding theorem of Epstein in \cite{Eps2}, we know $\phi$
is embedding.\stop


\subsection{Proof of Lemma \ref{Hu-Tal}}

\label{AppA}

In this subsection we prove Lemma \ref{Hu-Tal}. We start with \cite[Theorem 2.1]{Tal-1} as follows:

\begin{theorem}(\cite[Theorem 2.1]{Tal-1}) \label{Taliaferro 2006} Let $v(x,y)$ be a $C^{2}$
positive solution of 
\[
0\leq-\Delta v\leq Ce^{2v}
\]
 in a punctured neighborhood of the origin in $\mathbb{R}^{2}$ for a constant $C$. Then
either $v$ has $C^{1}$ extension to the origin or 
\[
\lim_{|x|\rightarrow0^{+}}\frac{v(x,y)}{\log(1/\sqrt{x^{2}+y^{2}})}=m_1
\]
for some finite positive number $m_1.$
\end{theorem}

\begin{remark}By considering $v(x,y)-\inf v(x,y)$, one can easily extend
the above theorem to the case that $v(x,y)$ is just bounded from below.
\end{remark}

To apply Theorem \ref{Taliaferro 2006} we first take an inversion. 
Let $\tilde{z}=\frac{z}{|z|^{2}}$
be the inversion map. Then 
\[
|dz|^{2}=\frac{1}{|\tilde{z}|^{4}}|d\tilde{z}|^{2} \ 
\text{and} \
g=e^{2u}|dz|^{2}=e^{2(u-2\log|\tilde{z}|)}|d\tilde{z}|^{2} = e^{2v}|d\tilde z|^2
\]
where 
\begin{equation}\label{Equ:v}
v(\tilde{z})=u(\frac{\tilde{z}}{|\tilde{z}|^{2}})-2\log|\tilde{z}|.
\end{equation}
We then have 
\begin{equation}\label{tilde-Gauss}
-\tilde{\Delta}v=K_{g}(\frac{\tilde{z}}{|\tilde{z}|^{2}})e^{2v} = \tilde K_g e^{2v} \quad\text{in \ensuremath{\mathbb{R}^2\backslash\{0\}}}.
\end{equation}

It is clear that, in order to apply Theorem \ref{Taliaferro 2006}, we need to obtain a lower bound first for the conformal factor $v$. 
To this purpose we first observe that $e^{-v}$ is a subharmonic function on $(\mathbb{R}^2, e^{2v}|d\tilde z|^2)$, that is,
$$
\Delta_{g} e^{-v} = e^{-v}|\nabla_g v|^2  - e^{-v}\Delta_g v = e^{-v} (|\nabla_g v|^2 + \tilde K_g) \geq 0.
$$
To obtain the lower bound, we recall \cite[Theorem 1.2]{Li and Schoen}. To state their theorem 
we consider a Riemannian manifold $M$, $x_0\in M$, and a radius $r$ such that, if $M$ has no boundary, 
$r$ is less than half of the diameter of $M$; if $\partial M\neq \emptyset$, $r < \frac 15 \text{dist}(x_0, \partial M)$. 

\begin{theorem} (\cite[Theorem 1.2]{Li and Schoen}) \label{Li-Schoen} Suppose that $M^n$ is a Riemannian manifold with $Ric \geq - (n-1)k$. Let $x_0\in M$ and $r$ given as above. Then for a nonnegative subharmonic function $v$ we have, for a constant $C$ depending only on the dimension and any $\tau\in (0, \frac 12)$,
\begin{equation}\label{Equ:Li-Schoen}
\sup_{B_{(1-\tau)r}(x_0)} v^2 \leq \tau^{-C(1 + \sqrt{k}r)}\frac 1{\text{vol}(B_r(x_0))} \int_{B_r(x_0)} v^2dvol.
\end{equation}
\end{theorem}

We therefore have, for the conformal factor $v$ in \eqref{tilde-Gauss},
\begin{equation}\label{Equ:e^{-v}}
\aligned
\sup_{B_{(1-\tau)r}(x_0)} e^{-2v} & \leq \tau^{-C}\frac 1{\text{vol}_{g}(B_r(x_0))} \int_{B_r(x_0)} e^{-2v}
dvol_{g}\\
& \leq \tau^{-C}\frac {\text{vol}_{|d\tilde z|^2} (B_r(x_0))} {\text{vol}_{g}(B_r(x_0))}.
\endaligned 
\end{equation}

Fortunately, we have a non-collapsing result in dimension 2 from \cite[Theorem A]{Croke and Karcher} as follows:

\begin{theorem} (\cite[Theorem A]{Croke and Karcher}) \label{Croke-Karcher} If $(M^2, \ g)$ is complete and nonnegatively
curved, then there exists a constant $C(M)$ such that, for $r\leq 1$,  
\begin{equation}\label{Equ:croke-karcher}
\text{vol}_g (B_r(x)) \geq C(M) r^2.
\end{equation}
\end{theorem}

Thus, the fact that the conformal factor $v$ is bounded from below follows from \eqref{Equ:e^{-v}}, \eqref{Equ:croke-karcher}, and the fact that
$\text{vol}_{|d\tilde z|^2}(B_r(x_0))$ is bounded. In fact, in this way we may conclude that $v(x_0) \to \infty$ as $\tilde z(x_0)\to 0$.
\\

Now we are ready to finish the proof of Lemma \ref{Hu-Tal}.
\\

\noindent{\it Proof of Lemma \ref{Hu-Tal}} \quad According to Theorem \ref{Taliaferro 2006}, we get
\[
v(\tilde{z})=m_{1}\log\frac{1}{|\tilde{z}|}+o(\log\frac{1}{|\tilde{z}|}) \ \text{as $\tilde z\to 0$}
\]
for some constant $m_1 > 0$. Next we claim that $m_{1}\geq 1$ since the metric $g = e^{2v}|d\tilde z|^2$
is complete and noncompact at the origin.  
\\

Assume otherwise $m_1 <1$.  Then let $m_2 \in (m_{1}, 1)$ and $r_{s}$ be sufficiently small so that 
\[
v < m_{2}\log\frac{1}{|\tilde z|}\,\,\text{for all }\,\,0<|\tilde z|<r_{s},
\]
which implies
\[
\exp(v)<|\tilde z|^{-m_{2}}\,\,\text{for all }\,\,0<|\tilde z|<r_{s}
\]
and 
\[
\int_{0}^{r_{s}}\exp(v(t,0))dt<\int_{0}^{r_{s}}t^{-m_{2}}dt < \infty.
\]
This contradicts the assumption that the metric $g = e^{2v}|d\tilde z|^2$ is complete and noncompact at the origin.
\\

Therefore, from \eqref{Equ:v}, we have 
\[
u(z)=(2-m_{1})\log\frac{1}{|z|}+o(\log\frac{1}{|z|})\,\,\text{as}\,\,|z|\to\infty,
\]
where $m=2-m_{1}\leq1$. 
\\

To see $m>0$ when $g$ is nonnegatively curved and nonflat, we recall
\[
-\Delta u=K_ge^{2u}\geq0\quad\text{in \ensuremath{\r^{2}}}.
\]
Taking an approach similar to that in the proof of \cite[Lemma 3]{Tal},
for $0<r_{2}<r_{1}$, we have that 
\begin{equation}
r_{2}\bar{u}'(r_{2})=r_{1}\bar{u}'(r_{1})+\frac{1}{2\pi}\int_{r_{2}<|z|<r_{1}}K_ge^{2u},\label{average again}
\end{equation}
where 
\[
\bar{u}(r)=\frac{1}{2\pi}\int_{0}^{2\pi}u(r\cos\theta,r\sin\theta)d\theta.
\]
Then 
\[
|\bar{u}'(r)|\leq\frac{1}{2\pi}\int_{0}^{2\pi}|\nabla u(r\cos\theta,r\sin\theta)|d\theta
\]
and therefore 
\[
\lim_{r_{2}\to0^{+}}r_{2}\bar{u}'(r_{2})=0.
\]
Plugging this back into \eqref{average again}, we have that 
\[
r_{1}\bar{u}'(r_{1})=-\frac{1}{2\pi}\int_{|z|<r_{1}}K_ge^{2u}.
\]
Now, from $u=m\log\frac{1}{|z|}+o(\log|z|)$ as $|z|\to\infty$, it
follows that 
\[
\lim_{r_{1}\to\infty}r_{1}\bar{u}'(r_{1})=-m=-\int_{\r^{2}}K_ge^{2u}<0,
\]
as $K_{g}\geq0$ and is not identically $0$.\stop


{\footnotesize{}{}}{\footnotesize \par}

\end{document}